\documentclass[letterpaper, 10 pt, conference]{ieeeconf}  

\IEEEoverridecommandlockouts                              
\usepackage{color}
\usepackage{graphicx}
\usepackage{epsfig} 
\usepackage{amsmath}
\usepackage{flushend}

\usepackage{caption}
\usepackage{subcaption}
\usepackage{cite}

\title{\LARGE \bf
Stochastic MPC Design for a Two-Component Granulation Process
}

\author{Negar Hashemian and Antonios Armaou$^{1}$
\thanks{$ ^1$Corresponding author, email: {\tt\small armaou@engr.psu.edu}
Financial support from the National Science Foundation, CBET Award 12-64902 is gratefully acknowledged. }
\thanks{N. Hashemian and A. Armaou are with the Department
of Chemical Engineering, The Pennsylvania State
        University, University Park, PA 16802.}%
}

\begin{document}

\maketitle
\thispagestyle{empty}
\pagestyle{empty}


\begin{abstract}
We address the issue of control of a stochastic two-component granulation process in pharmaceutical applications through using Stochastic Model Predictive Control (SMPC) and model reduction to obtain the desired particle distribution. We first use the method of moments to reduce the governing integro-differential equation down to a nonlinear ordinary differential equation (ODE). This reduced-order model is employed in the SMPC formulation. The probabilistic constraints in this formulation keep the variance of particles' drug concentration in an admissible range. To solve the resulting stochastic optimization problem, we first employ polynomial chaos expansion to obtain the Probability Distribution Function (PDF) of the future state variables using the uncertain variables' distributions. As a result, the original stochastic optimization problem for a particulate system is converted to a deterministic dynamic optimization. This approximation lessens the computation burden of the controller and makes its real time application possible.
\end{abstract}

\section{INTRODUCTION}
There are many systems in different fields which consist of particle populations such as crystallization, polymerization, granulations and viral infections. The particle distribution in these systems are defined as a multi variable function of particle properties e.g. type, size and/or composition. Mostly, the governing equation for these dispersed systems includes a population balance resulting in an integro-differential equation that involves both integrals and derivatives of the unknown particle distribution function. This paper studies one of these particulate processes that has been enhanced for application in the pharmaceutical industry; two-component high shear granulation. In this process, the granules are stuck together and form bigger particles through use of inactive binder droplets called excipient. In an ideal granulation process the composition and size of produced granules is the same, however, in reality the particle size and composition are distributed over a range. The objective in this paper is shaping the particles distribution based on the desired characteristics and constraints in a stochastic environment. 

The coagulation rate of particles in this process is determined by a weighting function which appears in the integrals called {\it coagulation kernel}. More specifically, in this process, the kernel is a function of particles size. Solving this equation for specific kernels is discussed in the literature frequently \cite{bicomponent_sol1}, however, for a general nonlinear kernel, there is no analytical solution to the population balance equation. An alternative method proposed by Matsoukas {\it et al.} is the constant-Number Monte Carlo (cNMC) algorithm \cite{Mat}, in which there should always be a constant number of particles in the \textquotedblleft simulation box'' with a varying volume during the evolution. This method has less computational cost than the traditional Monte Carlo algorithm which studies a finite number of particles in a fixed-volume simulation box \cite{Smith}. However, still this technique is too slow and not useful in the model predictive control formulation. As a result, to describe the bulk statistics of the process, we employ method of moments which is a powerful technique in derivation of deterministic models.  In previous works, the authors employed Taylor and Laguerre polynomial expansions to derive a closed finite dimensional ODE system that models a two-component coagulation process \cite{Negar1,Negar2}. This approach results in a tractable ODE model of the process which is in agreement with the Monte Carlo simulation results of the process. Also, this reduced order model is used for online estimation of the internal dynamics of the process. In this manuscript, we march on to design a Stochastic Model Predictive (SMPC) controller for this type of process.

Receding horizon approach is a powerful technique used in different control and estimation applications \cite{fang1,fang2,hashemian1,hashemian2,hashemian3,reza}. However, the underlying system is exposed to stochastic parameters at the feed flow. As a result, the control formulation should consider these uncertainties in its structure. In the literature, there are two approaches to accomplish this issue; Robust model predictive control \cite{RobustMPC,RobustMPC2,RobustMPC3} and stochastic model predictive control. Robust model predictive control is a more conservative method which considers the worst scenario in the optimization procedure. However, this method, similar to other robust controls \cite{Farshid1,Vahid1,Vahid2}, deteriorates the overall controller's performance and also is applicable only for systems with bounded uncertainties. The alternative method, SMPC, considers soft constraints which limit the risk of violation by a probabilistic inequality \cite{SMPC1,SMPC2,SMPC3}. This manuscript employs the later approach to control the granulation process. More specifically, in the granulation process, it is important that the active ingredient of the granule be in the admissible range. Also, the variance of average drug mass for different random scenarios should be minimized. We consider these factors in the cost function of the MPC structure. In addition, there are some soft constraints on the system to prevent the composition variance of particles to violate the admissible range. However, SMPC results in an optimization programming that is hard to solve in general. 

In the literature there are several efforts of sampling-based techniques which generally are computationally expensive and restricted to convex problems. To tackle this issue, chance constrained MPC is addressed in \cite{ashkan1} where  optimal control  input  for  given stochastic  dynamical  system  is obtained to minimize  a  given  cost-function  subject  to  probabilistic constraints,  over  a  finite  horizon. Building on the theory of measures and moments,  a  sequence  of  finite  semidefinite  programmings are   provided,   whose   solution   is   shown   to   converge   to the  optimal  solution  of  the  original  problem \cite{ashkan1,ashkan2,ashkan3}.  Mesbah {\it et al.} employ the generalized Polynomial Chaos (PC) theory to convert an SMPC formulation to a deterministic one~\cite{Mesbah}. PC expansion is a probabilistic method which projects the model's output in terms of orthogonal basis functions of random inputs. This stochastic method maps the future state variables from the uncertainty parameters utilizing orthogonal polynomials. Then, the approximate of state variables' statistical moments are available using the derived polynomial coefficients. This manuscript employs the same approach to control the granulation process.  
 
In this article, first we represent the population balance equation for a two component coagulation process and briefly explain how to derive a reduced order model from the original one. Next, a proper formulation for the SMPC with probabilistic constraints is proposed to control the granulation process. Afterward the article derives statistical properties of predicted state variables using the distribution of uncertainty variables. The last section applies the proposed SMPC on a two-component granulation system in the presence of noise in the feed flow concentration and compares the simulation results of this approach with a typical Nonlinear Model Predictive Control (NMPC).

\section{Problem Formulation}
\label{sec:2}

In the high shear granulation, an inactive ingredient (excipient) is fed to the active powder particles (drug/solute) in a tank with blending tools. In two-component granulation modeling, particles are distinguished by an augmented vector ${r=(p,s)}$, where $p$ and $s$ denote the total mass and the drug mass of the particle, respectively. The two variable function $f(p,s)$ shows the population distribution over the mass and solute content of particles. Since always the solute mass is less or equal to the total mass of the particle, the population distribution function has a zero value when $s>p$. 

In this process, the probability function of collision between $particle 1$ and $particle 2$ and a larger particle formation is denoted by $k_{12}=k(r_1,r_2)$. This function, called coagulation kernel, determines the dynamic behavior of particles in the system. As a result, for a general kernel function, the rate of particles' density with size $r_1$ using mass balance is given by  \cite{Mat}:
\begin{equation}
\begin{split}
\frac{\partial {f(r_1)}}{\partial t} &=
\frac {1}{2} \int \limits_0 ^{r_1} k(r_1-r_2,r_2) f(r_1-r_2) \ f(r_2)dr_2 \\
&
-  \int \limits_0 ^\infty k(r_1,r_2) f(r_1) f(r_2)dr_2
\label {eq:mass_balance}
\end{split}
\end{equation}

\subsection{Process Control}
In this process control there are two important goals: (i) uniform particles (ii) reaching to the desired size and composition. As a result, we are interested to minimize the deviation of expected values of $s$ and $p$ from their corresponding desired values and the variance of drug content in granules. Additionally, adding a probabilistic constraint guarantees the particles' drug mass stays within admissible region. To formulate this objective function, first we define the mixed moments as follows:
\begin{equation}
M_{ij}=\int _0 ^\infty \!\!\! \int _0 ^p p^is^j  f(p,s)\ dp ds
\label{eq:moments}
\end{equation}
Therefore the average value of drug and total mass of particles in each evolution are $\bar{s}=\frac{M_{01}}{M_{00}}$ and $\bar{p}=\frac{M_{10}}{M_{00}}$, respectively. Also, the variance drug amount in granules is $\frac{M_{02}}{M_{00}}$. However, as mentioned before, this process is exposed to uncertain parameters at the feed flow. As a result, we are interested to account for the expected values of these variables in the cost function instead of the specific values for a predefined input signal:

\subsubsection*{Problem 1 (Stochastic MPC with probabilistic constraints)}
\begin{equation}
\begin{split}
\min_{u_f} \bigg ( [E(\frac{M_{01}}{M_{00}})&-\mathcal{S}]^2 +[E(\frac{M_{10}}{M_{00}})-\mathcal{P}]^2+ \sigma Var(\frac{M_{01}}{M_{00}}) \bigg )\\
\mbox{Subject to: } &\mbox{Coagulation process model in Eq. (\ref{eq:mass_balance})}\\
& Pr [p_1^* \leq \frac{M_{02}}{M_{00}} \leq p_2^*] \geq \epsilon
\end{split}
\label{eq:sto_costfunction}
\end{equation}
where $\mathcal{S}$ and $\mathcal{P}$ are the desired mass of drug content and the total mass of the particles, respectively, $\sigma$ is a positive weight factor and $p_1^*,p_2^*$ are the lower and upper bound of admissible range for the variance of particle size, respectively. Moreover, $\epsilon \in (0,1)$ denotes the lower bound of the desired joint probability that particles' size should satisfy under uncertainties and $u_f$ is the manipulated variable in the process.

Additionally, the deterministic form  of the probabilistic constraint is given by \cite{probabilisticconstraints}:
\begin{equation}
\begin{split}
\kappa &Var(\frac{M_{02}}{M_{00}})-E(\frac{M_{02}}{M_{00}})+p_1^* \leq 0 \\
\kappa &Var(\frac{M_{02}}{M_{00}})+E(\frac{M_{02}}{M_{00}})-p_2^* \leq 0
\end{split}
\end{equation}
where $\kappa=\sqrt{\epsilon/(1-\epsilon)}$. 

To solve the above dynamic optimization problem, we need to solve the integro-differential equation at every sampling time and for each candidate guess. However, this equation does not have an explicit analytical solution for a nonlinear kernel function. In the literature, there are two general approaches to obtain the solution of these dynamic particulate systems; Monte Carlo simulation methods and conversion to an ODE set. The discrete nature of Monte Carlo simulation is very helpful to study the dynamical behavior of the particles, however, this approach is computationally expensive and inappropriate for estimation and control applications. Between the second method group of classification, we employ the method of moments to reduce the order of Eq.~(\ref{eq:mass_balance}). This approach not only simplifies the original computations, but also presents an approximation of population's probabilistic moments used directly in the SMPC formulation.

\subsection*{Model Order Reduction}
The method of moments obtains a reduced order model which gives the rate of change of mixed moments, $M_{ij}$. This model can be obtained by multiplying Eq.~(\ref{eq:mass_balance}) by $p^i s^j$ and then integrating over the region used in the definition of moments in Eq. (\ref{eq:moments}). However, because of the nonlinearities in the granulation kernel function, there exists no analytical solution of the double integral in the resulting equations by applying method of moments. To approximate this double integral, in our previous works\cite{Negar1,Negar2}, we used Taylor and Laguerre polynomial expansions. This method approximates the mixed moments \cite{Negar1} or directly approximates the population distribution versus a finite set of mixed moments\cite{Negar2}. These approaches result in two different deterministic ODE equation sets for Brownian kernel and kernels obtained from kinetic theory of granular flow (KTGF) \cite{KTGF}. In this manuscript we employ these reduced order models in the MPC structure to reduce the computation required in optimization stages. In the next section, we provide a method to approximate the expected values and variance of these moments of the system by considering the distribution of the uncertainties or noise signals.

\section{State Variables Statistics}

The previous section defined an objective function with probabilistic constraints. Then, we developed equivalent algebraic constraints for the stochastic optimization problem. However, still the expected values and variance of the process' state variables, i.e. mixed moments are required to control. The numerical calculation of expected value and variance of state variables by sampling noise distribution at every iteration makes the online solution of the problem computationally infeasible. As a result, this section employs the method introduced in \cite{Mesbah} to obtain a deterministic model to describe the statistical dynamic of the original mixed moments. This method uses polynomial chaos approximation to map the uncertainties on the dynamic system. Therefore, we predict the future state variables' PDFs using the uncertainties' PDFs. 

The reduced order model derived using method of moments in the previous section can be presented in the following compact form:
\begin{equation}
\begin{array}{l}
\displaystyle \dot{x}=f(x)+(B+w)\mu \\
\displaystyle y=x_1+\nu
\end{array} 
\label{eq:ODE}
\end{equation}
where $x, \mu$ are the vectors of process and input's distribution moments, $x_1$ is the first moment denoted by $M_{00}$ earlier. In this notation, $w$ and $v$ represent the uncertainty about the concentration of the particles entering the system and noise measurement at the output. 

To obtain the expected value and variance of these state variables, we use $s$ samples drawn from the known PDF of $w$ and obtain the corresponding state variables, in case these scenarios happen. Then, the PC expansion with the following structure approximates the stochastic state variables:
\begin{equation}
\hat x_t(w) = \sum_{i=0}^\infty a_i {\boldsymbol \phi_{\alpha,i}}(w)
\end{equation}
where $\boldsymbol \phi_{\alpha,i}(w)=\Pi_{j=1}^n \phi_{\alpha_{j,i}}(w_j)$ denotes the multivariate polynomials and the corresponding coefficient is given by $a_i=\frac{E[\hat x_t (w) {\boldsymbol \phi_{\alpha,i}}(w)]}{E[{\boldsymbol \phi_{\alpha,i}}(w)^2]}$. Also, $w=\begin{bmatrix}
w_1 &w_2  & \dotsb & w_n
 \end{bmatrix}$ is the vector of uncertain system parameters and $\phi_{\alpha_{j,i}}$'s are univariate polynomials chosen with respect to the distribution of the corresponding $w_j$. For example, the preferred choices for Gaussian, uniform and Gamma random variables are Hermite, Legendre and Laguerre polynomial bases, respectively.

Let $m_i=\sum_{i=1}^n \alpha_{j,i}$ and the order of the polynomial expansion is denoted by $m \geq m_i$. Assuming the ordering of polynomials satisfies the inequality $m_i \leq m_{i+1}$, therefore, the total number of polynomials in the truncated expansion is given by\cite{fagiano2012nonlinear}:
$$L=\frac{(n+m)!}{n!m!}$$
and the truncated expansion is:
\begin{equation}
\hat x_t(w) = \sum_{i=0}^{L-1} a_i {\boldsymbol \phi_{\alpha,i}}(w)
\end{equation}

Using orthogonality properties of the polynomial, the first and second moments of the stochastic state variables are given by:
\begin{equation}
\begin{array}{l}
\displaystyle E[\hat x(w)]=a_0 \\
\displaystyle Var[\hat x(w)]=\sum_{i=1}^{L-1}a_i^2 E[\phi_{a,i}^2(w)]
\end{array} 
\end{equation}
The higher order moments are obtained in terms of the polynomial coefficient in \cite{highermoments}. Also, we numerically calculate $E[\phi_{a,i}^2(w)]$ offline.
 

\section{Results and Discussions}
\label{sec:6}
This section investigates the nonlinear SMPC in the continuous two-component granulation process. In this section, first we use the PC Hermite expansion to approximate the dynamic behavior of statistic moments. Also, to assess the accuracy of this approach, the PDF of future state variables are obtained as a benchmark using a large enough number of noise samples in the feed flow. In the second part, the SMP structure is employed to reach the desired distribution in presence of noise in the feed flow and its performance is investigated.

\subsection*{PC Hermite expansion of granulation process}
This part evaluates the PC approximation performance assuming there exists a Gaussian noise on the feed flow drug mass. The reduced order model from \cite{Negar1} is employed. In this model, $C_f$ is the number of particles in the input flow and $\alpha$ denotes the input/output flow rate, where both the quantities are normalized by the coagulation container's volume. Additionally, it is assumed the aggregation kernel is Brownian with constant coefficient $k_0$. In all simulations, $k_0=0.06$, $\alpha=0.5, C_f = 1$ and a feed flow is assumed with uniform particles of total mass $p=1$. The feed flow concentration, $C_f$, is corrupted by a Gaussian noise signal. During the evaluation of PC approximation, particles' drug content at the feed flow are considered constant equal to $0.1$ and the initial state variable is:
\begin{equation*}
\begin{split}
x_0 = & \begin{smallmatrix} [M_{00} & M_{10} & M_{01} & M_{11} & M_{20} & M_{02} & M_{12} & M_{21} & M_{22}\end{smallmatrix} ]\\
=& \begin{bmatrix}1.9 & 2.0 & 0.2 & 0.2 & 2.3 & 0.02 & 0.03 & 0.3 & 0.05\end{bmatrix}
\end{split}
\end{equation*}
The PC coefficients are estimated using the projection method, where the integral is solved using a Gauss-Hermite quadrature.

\begin{figure}[!htb]
	\vspace{1mm}
	\centering
	\begin{subfigure}[b]{.4\textwidth}
		\includegraphics[width=\textwidth]{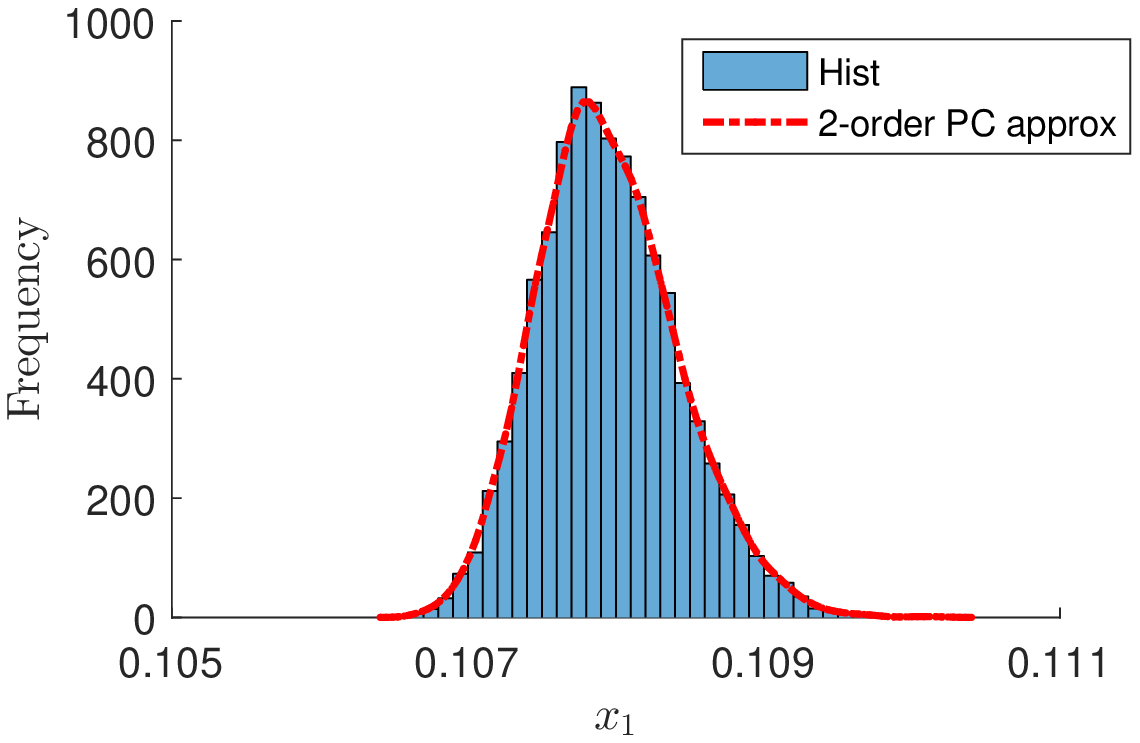}
		\caption{}
	\end{subfigure}
	\begin{subfigure}[b]{.4\textwidth}
		\includegraphics[width=\textwidth]{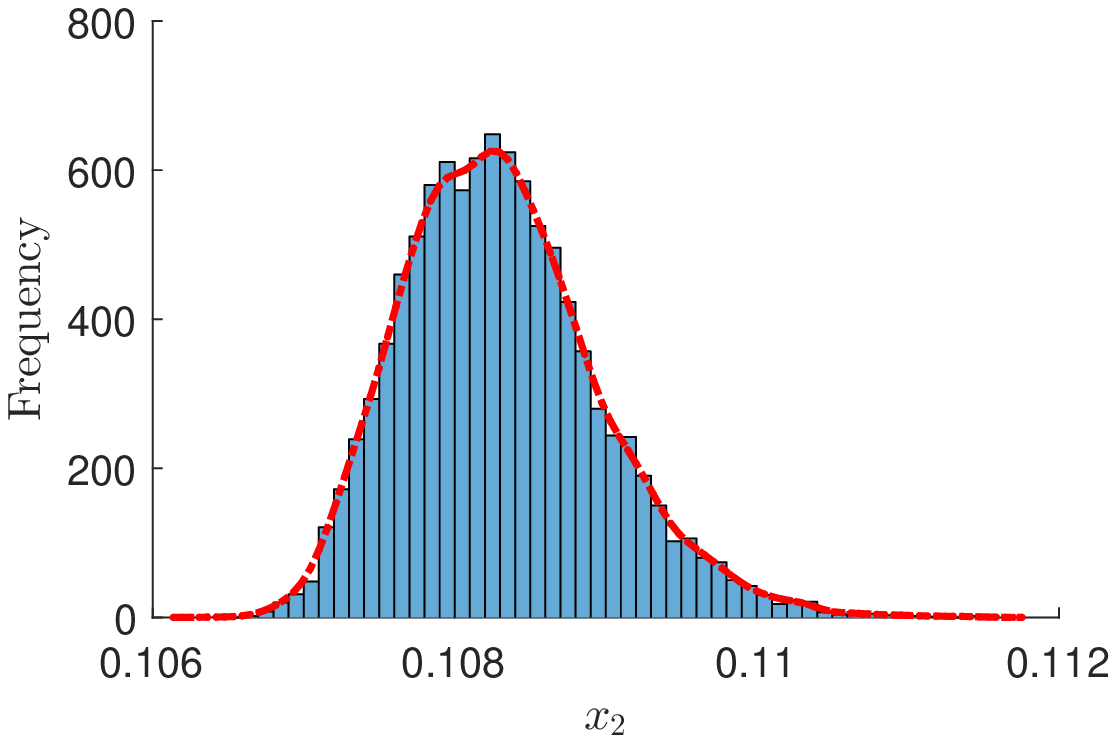}
		\caption{}
	\end{subfigure}
	
	\begin{subfigure}[b]{.4\textwidth}
		\includegraphics[width=\textwidth]{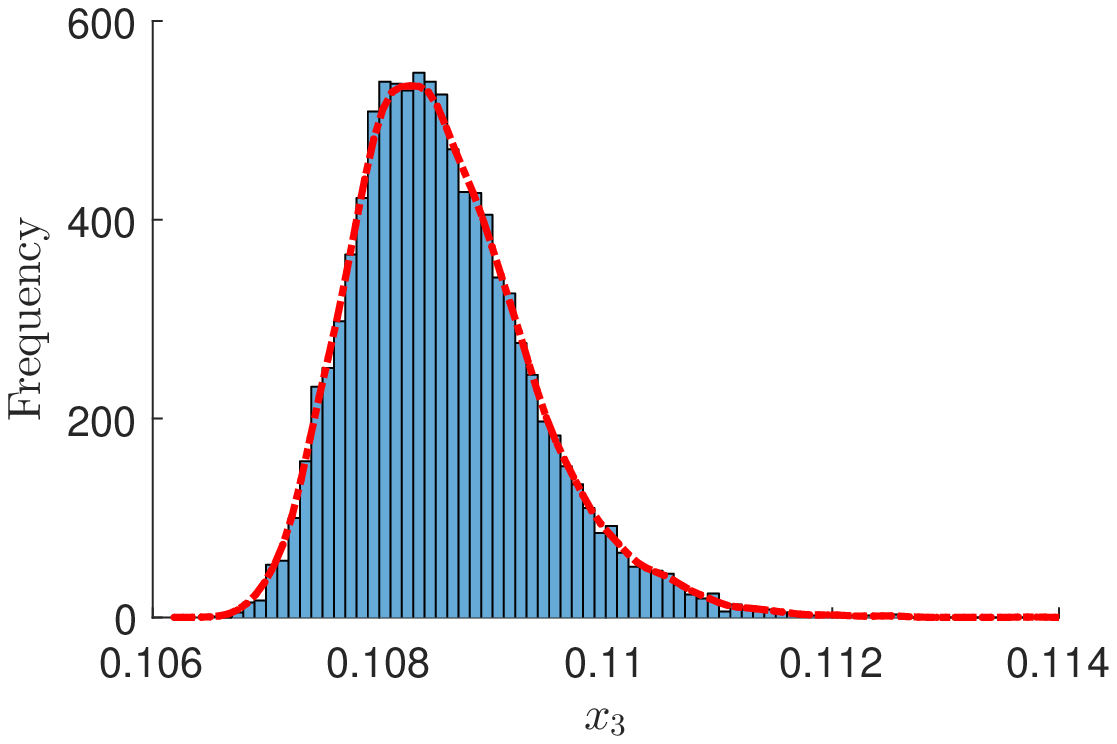}
		\caption{}
	\end{subfigure}
	\caption{Monte Carlo (MC) and Hermite PC solution of the average of the particles' drug content after (a) $\Delta T$, (b) 2$\Delta T$ and (c) 3$\Delta T$.}
	\vspace{-.75cm}
	\label{fig:PC}
\end{figure}

 The SMPC designed in the next part, has a prediction horizon $N=3$. As a result, the statistical variables are required to be approximated at the three future points  and the PC functions have the following structures:
\begin{equation}
\begin{split}
\hat{x}_{t+1}(w_1) &= a_0 H_i (w_1)+a_1 H_1 (w_1)+a_2 H_2 (w_1)\\
\hat{x}_{t+2}(w_1,w_2) &= \sum_{i=0}^{2} \sum_{j=0}^{2} a_{i,j} H_i (w_1) H_j (w_2)\\
\hat{x}_{t+3}(w_1,w_2,w_3) &= \sum_{i=0}^{2} \sum_{j=0}^{2} \sum_{k=0}^{2} a_{i,j,k} H_i (w_1) H_j (w_2) H_k (w_3)
\end{split}
\label{eq:HermitePlynomials}
\end{equation}
Fig. \ref{fig:PC} shows the distribution of $\frac{M_{10}}{M_{00}}$ after $\begin{bmatrix} \Delta T & 2\Delta T & 3\Delta T \end{bmatrix}$  denoted by $x_1, x_2$ and $x_3$ and $\Delta T = 0.25$. The results are obtained by the nonlinear model and PC expansion functions in Eq. (\ref{eq:HermitePlynomials}), whose coefficients are determined using 6 sample points. These points are the roots of the Hermite polynomials $H_i(w) \; (i = 0,1,, ... , 5)$. The PDFs obtained from the PC approximation are in the agreement with the histograms of dynamic variable with a random sample of $10,000$ points from a normal distribution for $\begin{bmatrix} w_1 & w_2 & w_3 \end{bmatrix}$.

\subsection*{Stochastic Control of the granulation process}

This section designs an SMPC structure using PC expansion method to control the continuous bi-component granulation process with fluctuations in the feed flow concentration. To compare the performance of this controller with a typical nonlinear model predictive control, Monte Carlo simulations of the closed-loop system are employed. These simulations are performed for 100 different input concentration trajectories sampled from a Gaussian distribution with a mean of $C_{fs}=1$ and a standard deviation of $0.1$. The input's particles are assumed to be uniform and the manipulated variable is the particles' drug content at the feed flow, $s_f$.  

\begin{figure}[!htb]
	\centering
	\begin{subfigure}[b]{0.50\textwidth}
		\includegraphics[width=\textwidth]{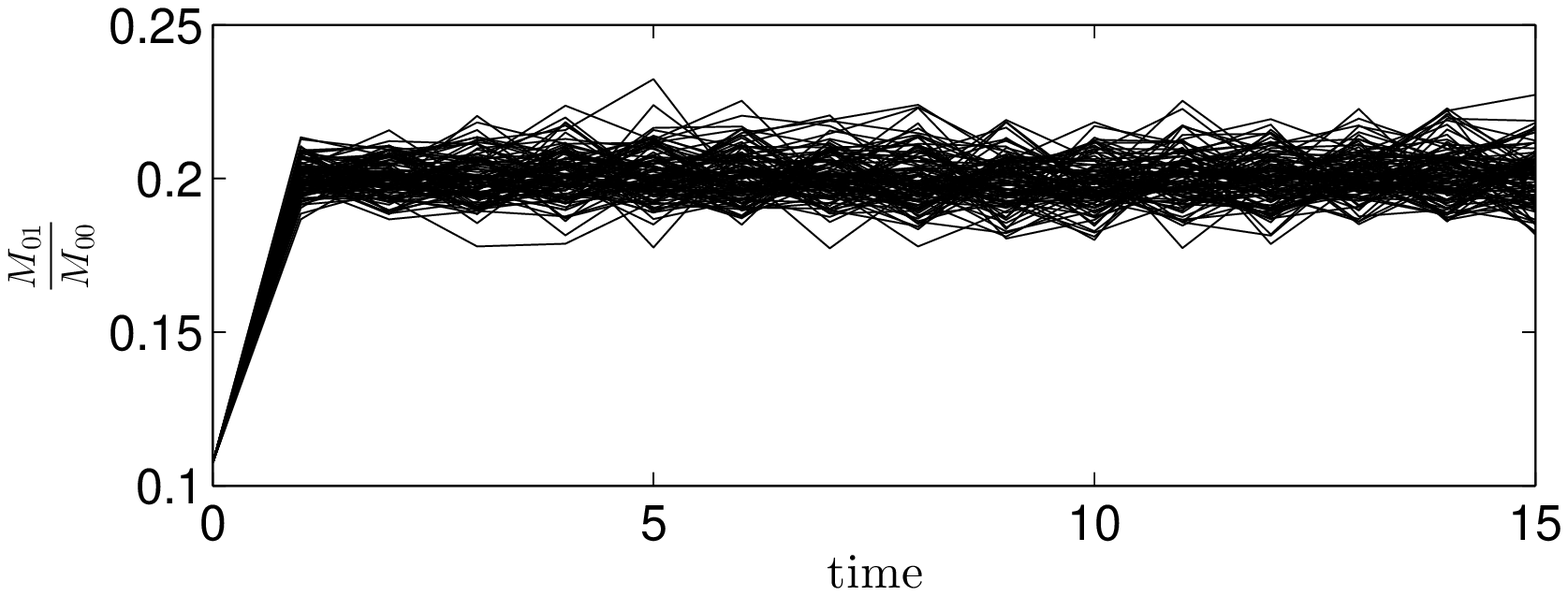}
		\caption{}
	\end{subfigure}
	\begin{subfigure}[b]{0.50\textwidth}
		\includegraphics[width=\textwidth]{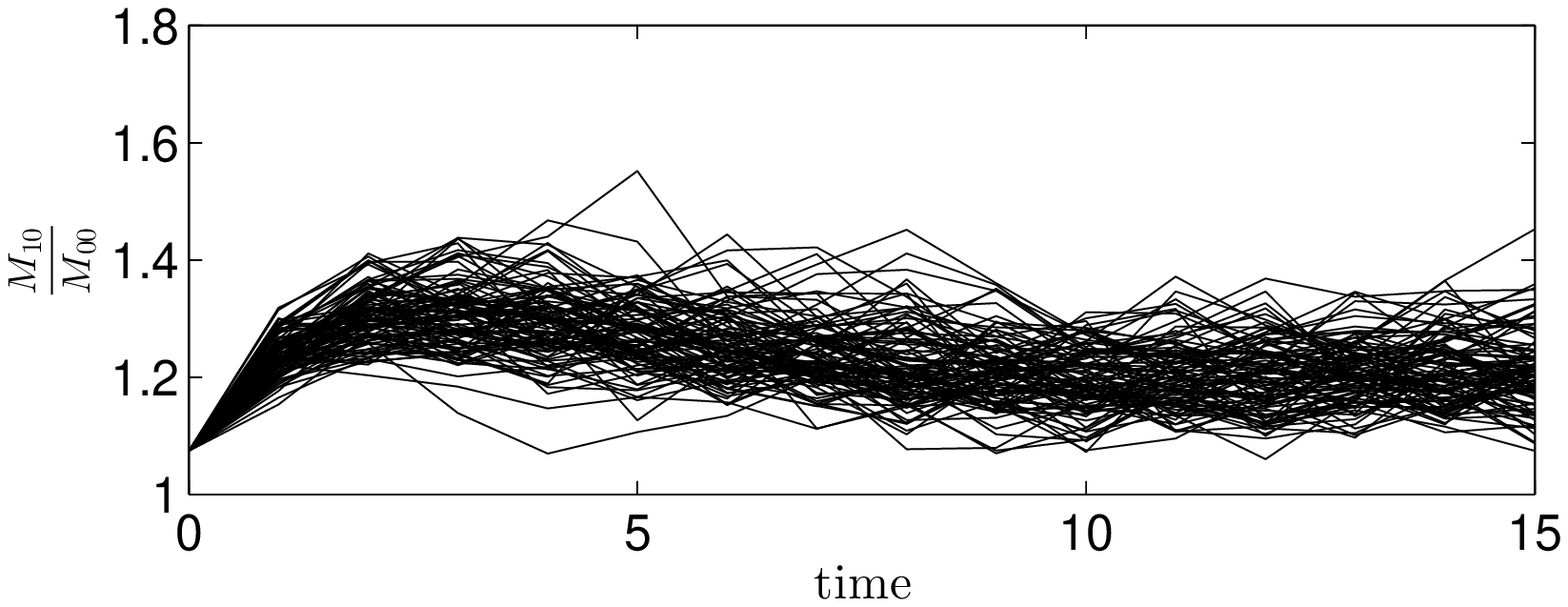}
		\caption{}
	\end{subfigure}
	
	\begin{subfigure}[b]{0.50\textwidth}
		\includegraphics[width=\textwidth]{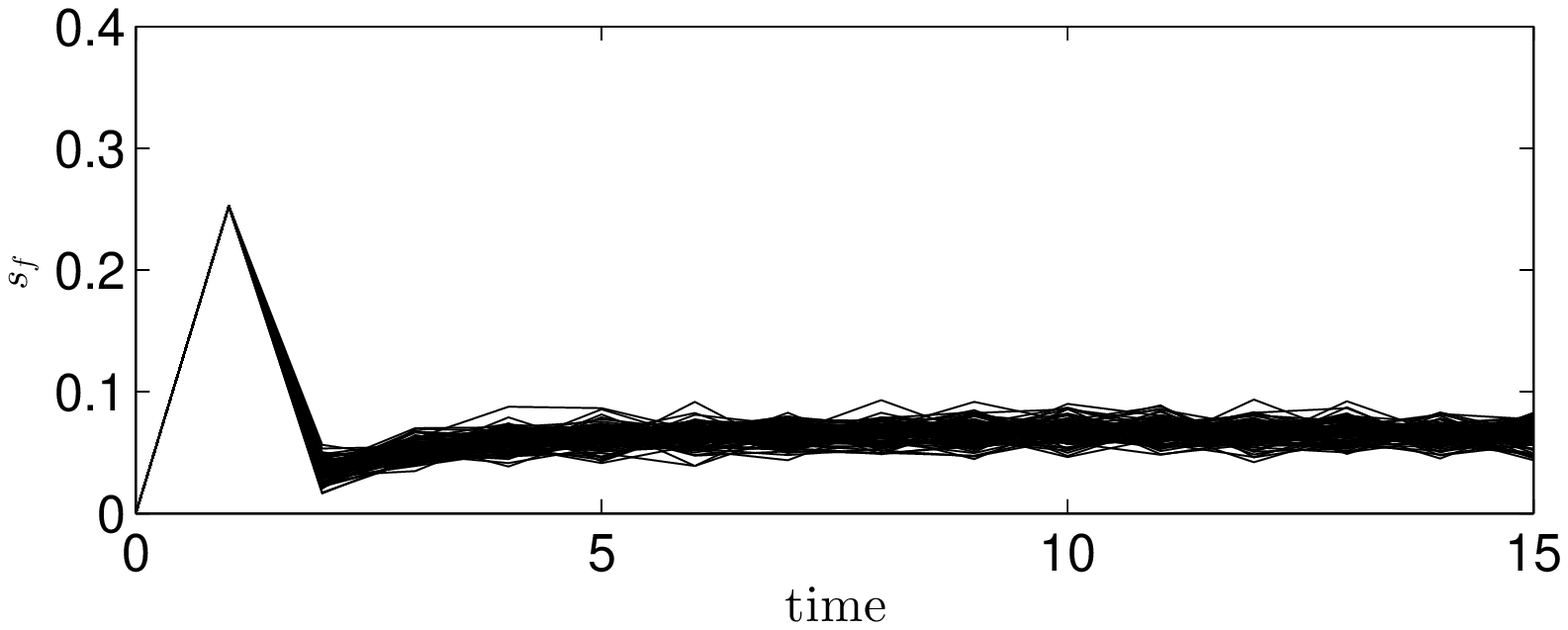}
		\caption{}
	\end{subfigure}
	\caption{ (a) Average particles' drug mass, (b) average particles' total mass and (c) input trajectories of the SMPC closed loop response.}
	\label{fig:SMPC}
	\vspace{-.6cm}
\end{figure}

\begin{figure}[!htb]
	\centering
	\begin{subfigure}[b]{0.45\textwidth}
		\includegraphics[width=\textwidth]{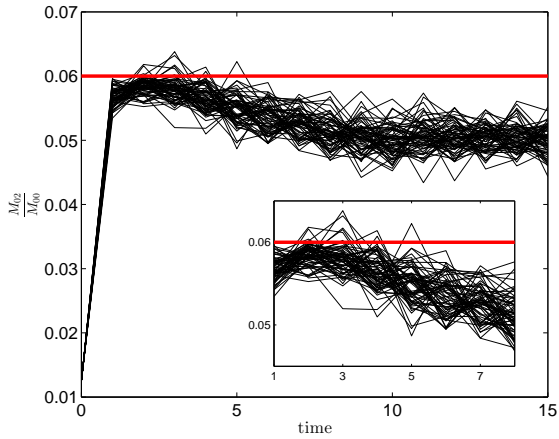}
	\end{subfigure}
	\caption{Drug mass variance trajectories for 100 simulations using SMPC design. The inset shows a magnification of $\frac{M_{02}}{M_{00}}$ to
		depict constraint violations.}
	\label{fig:cons}
	\vspace{-.6cm}
\end{figure} 

In all simulations, $\mathcal{S}=0.2, \mathcal{P}=1.2, p_1^*=0, p_2^*=0.06$ and the sampling time, $\Delta T=1$. Additionally, in the SMPC formulation, $\sigma = 100$ and $\epsilon = 0.85$. Figure \ref{fig:SMPC} shows the control input and the closed-loop system response. The controller is able to reject the disturbance in the flow rate and steers the system close to the desired particles properties at the output. To solve all the dynamic optimization problems, fmincon function in MATLAB with `interior-point' algorithm is utilized.
Moreover, Fig. \ref{fig:cons} shows the variance of drug mass in particles during the system evolution. The variable $\frac{M_{02}}{M_{00}}$ may violate its upper bound with a likelihood of $10\%$ in the presence of feed flow uncertainties. As shown in the inset figure, the constraint violation occurs only in the time between $t=1$ and $t=5$.

\begin{figure}[!htb]
	\vspace{1mm}
	\centering
	\begin{subfigure}[b]{.46\textwidth}
		\includegraphics[width=\textwidth]{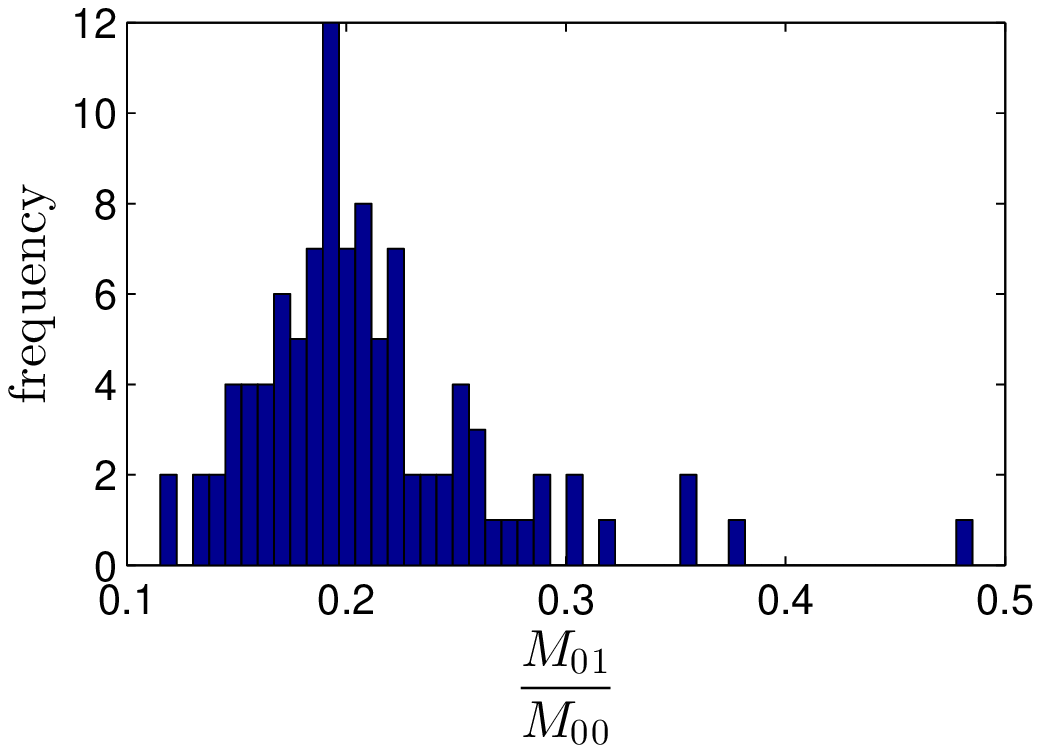}
		\caption{}
	\end{subfigure}
	\begin{subfigure}[b]{.46\textwidth}
		\includegraphics[width=\textwidth]{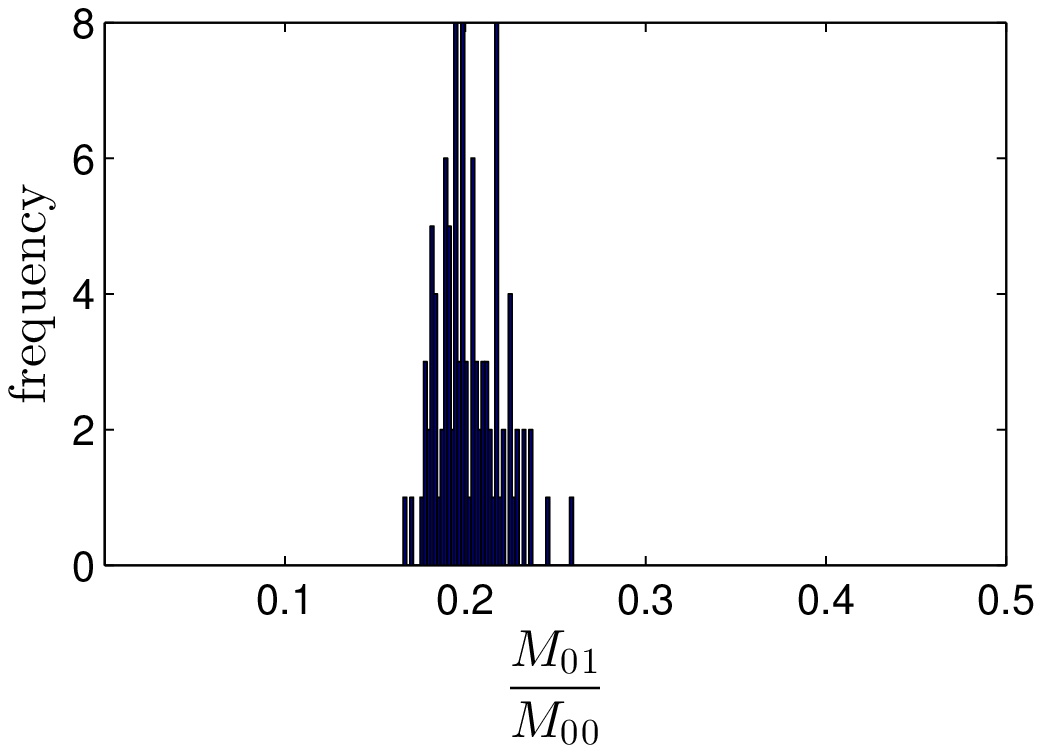}
\caption{}
	\end{subfigure}
\caption{Histograms of the predicted average of particles' drug content  at $t=15$ based on 100 simulations using (a) NMPC (b) SMPC design.\vspace{-.6cm}}
\label{fig:comparison1}
\end{figure}
	
	\begin{figure}[!htb]
		\vspace{1mm}
		\centering
	\begin{subfigure}[b]{.43\textwidth}
		\includegraphics[width=\textwidth]{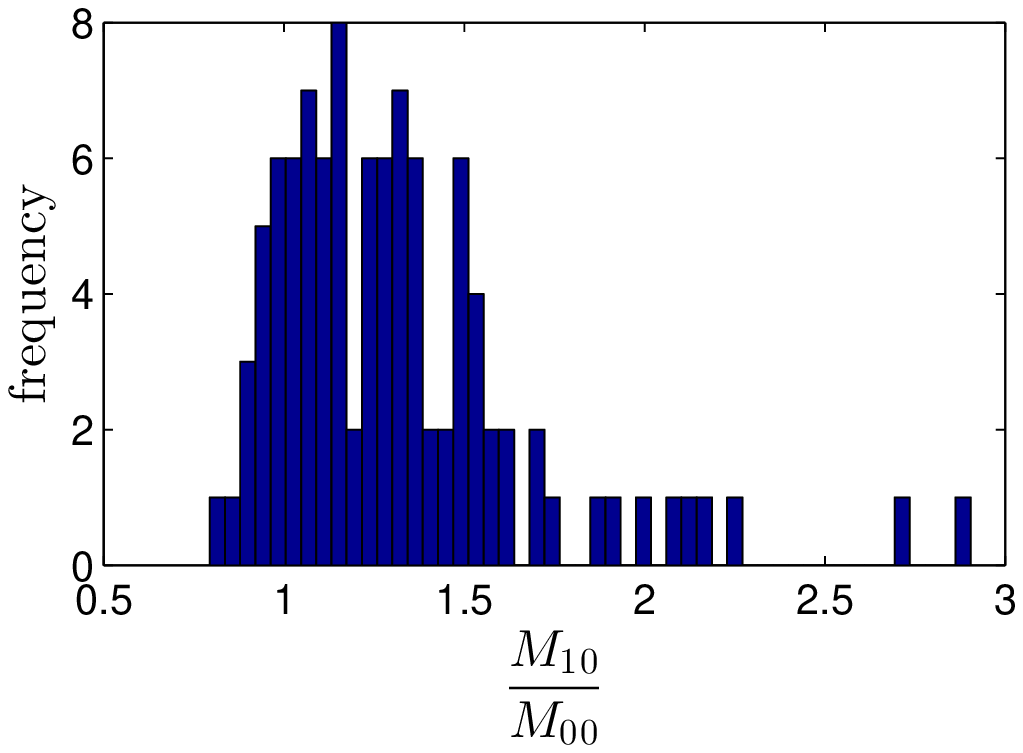}
		\caption{}
	\end{subfigure}
	\begin{subfigure}[b]{.43\textwidth}
		\includegraphics[width=\textwidth]{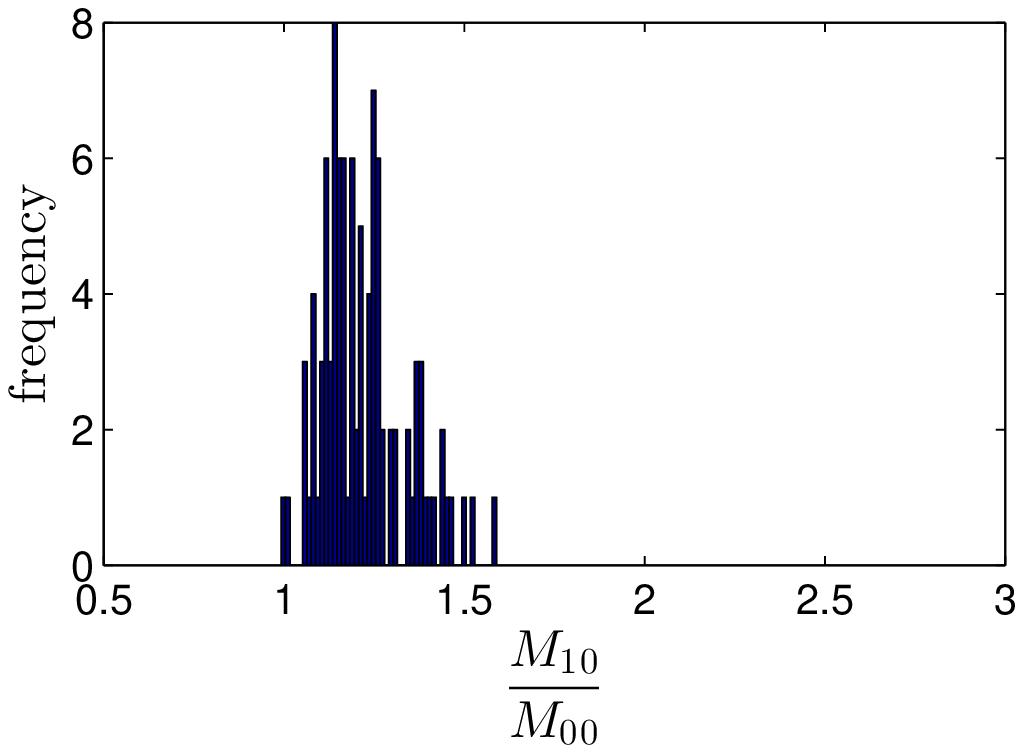}
		\caption{}
	\end{subfigure}
\caption{Histograms of the predicted average of particles mass at $t=15$ based on 100 simulations using (a) NMPC (b) SMPC design.}
\label{fig:comparison2}
\vspace{-.6cm}
\end{figure}

To investigate the performance of the SMPC design, an NMPC structure is formulated as follows:
\begin{equation}
\begin{split}
\min_{s_f} \bigg ( (\frac{M_{01}}{M_{00}}&-\mathcal{S})^2 +(\frac{M_{10}}{M_{00}}-\mathcal{P})^2 \bigg )\\
\mbox{Subject to: } &\mbox{reduced order model in Eq. (\ref{eq:ODE})}\\
& p_1^* \leq \frac{M_{02}}{M_{00}} \leq p_2^*
\end{split}
\label{eq:NMPC}
\end{equation}
 Fig. \ref{fig:comparison1} depicts the histograms of $\frac{M_{01}}{M_{00}}$ based on the Monte Carlo closed-loop simulations of the two control approaches. The distributions of the average particles' drug mass are with a mean of $0.2102$ and $0.2005$ and a variance of $3.2 \times 10^{-3}$ and $8 \times 10^{-5}$, in the NMPC and SMPC approaches, respectively. 

Fig. \ref{fig:comparison2} demonstrates the average particle mass at $t=15$, which is desired to be at $\mathcal{P}$. The results indicate that in
the NMPC approach $E[\frac{M_{10}}{M_{00}}]=1.3$ and $Var[\frac{M_{10}}{M_{00}}]=0.14$, while in the SMPC simulation $E[\frac{M_{10}}{M_{00}}]=1.212$ and $Var[\frac{M_{10}}{M_{00}}]=4.3 \times 10^{-3}$. Therefore, the simulation results of the two-component coagulation process shows SMPC approach with probabilistic constraints shapes the PDFs
of state variables properly and also satisfies the physical constraints on the system in the presence of stochastic uncertainties.

\section{Conclusions}
This article addresses the problem of stochastic nonlinear control of a two-component coagulation process. The distribution of particles in this process is described as a function of the particles' mass and composition. The particle balance for this system results in an integro-differential equation which does not have a closed form analytical solution in general. As a result, the model reduced order is exploited in the control formulation. 

A stochastic controller formulation is presented for the coagulation system to obtain the desired expected value of the drug and total mass with minimum variance of the drug content. Also, this stochastic dynamic optimization formulation keeps the violation probability of constraints in an admissible range. To simplify the optimization problem, the probabilistic inequalities are converted to some algebraic convex second order cone inequalities. Moreover, polynomial chaos expansions are exploited to predict the state variables distribution in presence of the feed flow's noise. The simulation results show the SMPC formulation shapes the probability distribution of system states, as well as guaranteeing the state constraints satisfaction in a stochastic environment.

\bibliographystyle{IEEEtran}
\bibliography{ref.Negar}

\end{document}